\documentclass[10pt,reqno]{amsart}

\usepackage[letterpaper,hmargin=1.25in,vmargin=1.22in]{geometry}
\usepackage{amssymb}
\usepackage[all]{xy}
\usepackage[hyphens]{url}
\usepackage{graphicx,xcolor,color}

\newtheorem{de}{Definition}[section]
\newtheorem{lem}[de]{Lemma}
\newtheorem{prop}[de]{Proposition}
\newtheorem{cor}[de]{Corollary}
\newtheorem{thm}[de]{Theorem}

\theoremstyle{remark}
\newtheorem{rem}[de]{Remark}
\newtheorem{ex}[de]{Example}

\def \N{\mathbb{N}}
\def \Z{\mathbb{Z}}
\def \Q{\mathbb{Q}}

\def \R{\mathbb{R}}

\newcommand{\Diff}{\mathfrak{D}\mathrm{iff}}
\newcommand{\Top}{\mathfrak{T}\mathrm{op}}

\newcommand{\cl}{\operatorname{cl}}
\newcommand{\Int}{\operatorname{int}}

\title{Topology on diffeological vector spaces}

\author{Enxin Wu}
\email{exwu@stu.edu.cn}
\address{Department of Mathematics, Shantou University, Shantou 515063, China}

\author{Zhongqiang Yang}
\email{zqyang@stu.edu.cn}
\address{School of Mathematics and Statistics,
Minnan Normal University, Zhangzhou 363000, China}

\thanks{Zhongqiang Yang is the corresponding author.}
\thanks{This work was supported by National Natural Science Foundation of China
	    (Nos. 11971141 and 11971287).}

\date{\today}

\begin{document}

\begin{abstract}
It is expected that the $D$-topology makes every diffeological vector space into a topological vector space. We show that it is the case for a large class of diffeological vector spaces via $k_\omega$-space theory, but not so in general. The paper also proposes the study of a class of almost topological vector spaces.
\end{abstract}

\makeatletter
\@namedef{subjclassname@2020}{%
  \textup{2020} Mathematics Subject Classification}
\makeatother


\maketitle

\section{Introduction}

A relatively new field named diffeology (introduced by J.-M. Souriau in the 1980's~\cite{S}) provides a uniform framework for smooth manifolds, infinite-dimensional spaces (such as function spaces and diffeomorphism groups) and singular spaces (eg. smooth manifolds with boundary or corners, orbifolds, irrational tori, etc). The principle is, only the classical topology and analysis on Euclidean open subspaces are admitted and are used as building blocks, and all other results have to be derived from strict rules. The definition of a diffeological space is very simple (see Definition~\ref{de:basics}\eqref{diffspace}). On such spaces, one can still do differential topology: tangent spaces and tangent bundles, differential forms, de Rham cohomology, fibre bundle theory, smooth homotopy groups, etc, which are natural generalizations of the corresponding concepts for smooth manifolds. Meanwhile, there is a unique non-trivial extreme topology (called the $D$-topology) on every diffeological space, first introduced in~\cite{I85}, which usually works as an important auxiliary tool.

A vector space object in the category of diffeological spaces is called a diffeological vector space. It appears in tangent bundles (\cite{CW16, CW17}), and more generally vector bundle theory (\cite{CW21,CWp}), and homological algebra (\cite{W}), etc. It is expected that the $D$-topology on every diffeological vector space makes it into a topological vector space, so that one has a bridge between diffeological vector space theory and topological vector space theory, and it was claimed so in \cite[Example~2.8]{W} without proof by the first author. This was implicitly based on the assumption that product and $D$-topology commute. However, this is not the case in general; see Section~\ref{s:notcommute}. Luckily this false claim was not used elsewhere in~\cite{W}. On the other hand, for many known examples, the $D$-topology does send the diffeological vector space to a topological vector space. In the current paper, we find a proper condition to make this happen (Theorem~\ref{thm:proj}), and also provide an example (Theorem~\ref{prop:surprise2}) showing that this is not always the case, even for a diffeological vector space which is both free and projective. (Notice that a free diffeological vector space is not necessarily projective.)

Technically speaking, the $D$-topology makes every diffeological vector space into an almost topological vector space. To be precise, scalar multiplication is always continuous, and addition is continuous in each variable, but not necessarily both. (Maybe it is worthwhile to study such almost topological vector spaces separately, which does not seem to exist in the existing literatures related to generalizations of topological vector spaces.) Finally the focus is on the confirmation of the sufficient conditions when product and $D$-topology commute, whence the appearance of $k_\omega$-space theory.

General diffeological vector space is not easy to deal with directly. The strategy goes as follows. \cite{W} showed that every diffeological vector space is a quotient of a projective diffeological vector space (i.e., enough projectives), and every projective diffeological vector space is a smooth direct summand of direct sum of free diffeological vector spaces $F(U)$ generated by open subsets $U$ of Euclidean spaces (i.e., the structural theorem of projective diffeological vector spaces). Via $k_\omega$-space theory (Section~\ref{s:komega}), we show that $F(U)$ with its $D$-topology is a topological vector space (Section~\ref{s:apply}), and then show that this can be carried over to countable direct sum and smooth direct summand, and hence the same result is valid for countably-generated projective diffeological vector spaces via the structural theorem (Section~\ref{s:projective}). Using a concrete example that product and $D$-topology does not commute (Section~\ref{s:notcommute}) together with the technique of checking $D$-closedness throughout Section~\ref{sss:general}, we provide a free diffeological vector space generated by a smooth Tychonoff space which is uncountably-generated and projective, such that it is not a topological vector space under the $D$-topology (Theorem~\ref{prop:surprise2}). Furthermore, from this, we get some further related results.

To avoid mistakes or misuse of abstract results, this paper is self-contained except quoting several well-known results where refernces are provided in detail. In particular, we provide a detailed and concise introduction to the basic theory of $k_\omega$-spaces needed later in Section~\ref{s:komega}, since this material is not so well-known outside general topology. So there are some overlap with the existing literature on free topological vector spaces when dealing with the $D$-topology on $F(U)$, while we start from a quite different point of view. As a byproduct of our treatment, the paper contains an alternative characterization of the topology of a free topological vector space in good cases (Section~\ref{s:comparison}).

\section{Preliminaries}\label{s:background}

In this section, we collect the basics of diffeological spaces and diffeological vector spaces. \cite{I13} is currently the standard textbook on the subject, and~\cite[Section~2]{CSW} provides a concise introduction to the details of the constructions. See~\cite{CSW} for more on the $D$-topology, and see~\cite{W,CW19,Wp} for more on diffeological vector spaces.

\begin{de}\label{de:basics}\
\begin{enumerate}
\item\label{diffspace} A \textbf{diffeological space} $(X,\mathcal{D}_X)$ is a set $X$ together with a collection $\mathcal{D}_X$ of maps $U \to X$ (called \textbf{plots}), where $U$ is open in $\R^n$ and $n \in \N$, such that the following conditions hold:
\begin{enumerate}
\item Every constant map is in $\mathcal{D}_X$.
\item If $U \to X$ is in $\mathcal{D}_X$ and $V \to U$ is a smooth map between open subsets of Euclidean spaces, then the composite $V \to U \to X$ is also in $\mathcal{D}_X$.
\item $U \to X$ is in $\mathcal{D}_X$ if there exists an open cover $U = \bigcup_i U_i$ of $U$ such that each restriction $U_i \to X$ is in $\mathcal{D}_X$.
\end{enumerate}
Whenever $\mathcal{D}_X$ is understood, we use $X$ to denote the diffeological space $(X,\mathcal{D}_X)$.
\item A \textbf{generating set} of plots for a diffeological space $X$ is a subset $\mathcal{G}$ of $\mathcal{D}_X$ such that for every $p:U \to X$ in $\mathcal{D}_X$ and every $u \in U$, there exists an open neighborhood $U'$ of $u$ in $U$, so that either $p|_{U'}$ is constant, or there exist $V \to X$ in $\mathcal{G}$ and a smooth map $U' \to V$ with $p|_{U'}$ being the composite $U' \to V \to X$.
\item\label{smooth} A function $f:X \to Y$ between diffeological spaces is called \textbf{smooth} if $f(\mathcal{D}_X) \subseteq \mathcal{D}_Y$.
\item The \textbf{$D$-topology} on a diffeological space $X$ is the finest topology on the underlying set $X$ such that every plot (with domain equipped with the Euclidean topology) is continuous. We use $D(X)$
    to denote the topological space of $X$ with the $D$-topology.
\end{enumerate}
\end{de}

Diffeological spaces together with smooth maps form a category, denoted by $\Diff$. It has the following convenient properties:

\begin{thm}\label{thm:basics}\
\begin{enumerate}
\item The category of smooth manifolds\footnote{As the usual convention, a smooth manifold is assumed to be finite-dimensional, Hausdorff, second-countable and without boundary.} and smooth maps is a full subcategory of $\Diff$. (In particular, smoothness in the sense of Definition~\ref{de:basics}\eqref{smooth} between smooth manifolds coincides with the smoothness in the usual sense. Without specification, every smooth manifold is equipped with this standard diffeology.)
\item $\Diff$ is bicomplete and cartesian closed.
\end{enumerate}
\end{thm}

In this way, we can define subspace and quotient space of a diffeological space, together with arbitrary product and function space of diffeological spaces. A morphism in $\Diff$ that is equivalent to a quotient map is called a \textbf{subduction}. It is well-known that the category of smooth manifolds is not bicomplete nor cartesian closed; see~\cite[Section~2]{GW} and references therein for a detailed discussion. So the category of diffeological spaces can be viewed as a remedy with respect to set-theoretical operations.

Here are the basic properties of the $D$-topology that we will use frequently in this paper:

\begin{prop}\
\begin{enumerate}
\item The $D$-topology on a smooth manifold coincides with the usual topology.
\item The $D$-topology gives rise to a left adjoint $D:\Diff \to \Top$, where $\Top$ denotes the category of all topological spaces and continuous maps. (In particular, the $D$-topology makes every smooth map between diffeological spaces to be continuous, and the $D$-topology functor commutes with colimits.)
\item The $D$-topology behaves not so nicely with limits. For example, if $A$ is a subset of a diffeological space $X$, then the $D$-topology of the diffeological subspace $A$ contains the topology of the corresponding subspace of the space $D(X)$, but they are different in general. Also, the $D$-topology of a finite product of diffeological spaces contains the product topology of the $D$-topology of each space, but they are equal only under strong conditions. See~\cite{CSW} for more details. In the current paper, we have more related results in Sections~\ref{s:apply}, \ref{s:projective}, and~\ref{s:general} for diffeological vector spaces.
\item Let $\mathcal{G}$ be a generating set of plots for a diffeological space $X$. Then $A \subseteq X$ is $D$-open (resp. $D$-closed) if and only if $p^{-1}(A)$ is open (resp. closed) for every $p \in \mathcal{G}$.
\end{enumerate}
\end{prop}

In this paper, every vector space and every linear map is over the field $\R$ of real numbers.

\begin{de}
A \textbf{diffeological vector space} $V$ is both a diffeological space and a vector space such that the addition $V \times V \to V$ and the scalar multiplication $\R \times V \to V$ are both smooth in the sense of Definition~\ref{de:basics}\eqref{smooth} and Theorem~\ref{thm:basics}.
\end{de}

In the current paper, we try to answer the question that under what condition does the $D$-topology make a diffeological vector space $V$ into a topological vector space. The obstacle comes from the fact that the natural continuous map $D(X \times Y) \to D(X) \times D(Y)$ induced by the projections $X \times Y \to X$ and $X \times Y \to Y$ may not be a homeomorphism for general diffeological spaces $X$ and $Y$; see Section~\ref{s:notcommute}. That is, generally speaking, the  $D$-topology of the product $X\times Y$ of two diffeological spaces $X$ and $Y$ is finer than the product topology of the
$D$-topology of the diffeological space $X$ and the
$D$-topology of the diffeological space $Y$. Note that~\cite[Lemma~4.1]{CSW} tells us $D(\R \times V) = D(\R) \times D(V)$, and apparently every translation $D(V) \to D(V)$ is continuous (i.e., addition $D(V) \times D(V) \to D(V)$ is continuous in each variable). So the key is to check whether the addition $D(V) \times D(V) \to D(V)$ is continuous, although the addition $D(V \times V) \to D(V)$ induced by $+:V \times V \to V$ is.\footnote{Notice a key difference that product of two quotient maps in $\Diff$ is again a quotient map, while this is not necessarily the case in $\Top$.}

\section{Product does not commute with $D$-topology}\label{s:notcommute}

The following example is adjusted from~\cite[Example in 5.1]{D}, which shows that $D(X \times Y) \neq D(X) \times D(Y)$ for general diffeological spaces $X$ and $Y$. The extensions in the successive remark will be used later.

\begin{ex}\label{ex:notcommute}
Let $X$ be the wedge (with the quotient diffeology) at $0$ of copies of $\R$ indexed over all infinite sequences $s = (s_1,s_2,\ldots)$ of positive integers, and let $Y$ be the wedge at $0$ of copies of $\R$ indexed over all positive integers.  We claim that $D(X \times Y) \neq D(X) \times D(Y)$.

Here is the proof. Write $p_{sj} := (1/s_j, 1/s_j)$ in the set $\R_s \times \R_j \subset X \times Y$, and let $P = \{p_{sj}\}_{s,j}$. Then
\begin{equation}\label{eq:closed}
P \text{ is closed in } D(X \times Y).
\end{equation}
This is because we have quotient maps $\coprod_s \R_s \to X$ and $\coprod_j \R_j \to Y$, which induce a quotient map $\pi:(\coprod_s \R_s) \times (\coprod_j \R_j) \to X \times Y$ in $\Diff$. By cartesian closedness of $\Diff$, we get $(\coprod_s \R_s) \times (\coprod_j \R_j) = \coprod_{s,j} \R_s \times \R_j$. Using the left adjointness of the $D$-topology functor, we know that $\pi^{-1}(X \times Y \setminus P)$ is open in $\coprod_{s,j} \R_s \times \R_j$, and hence $X \times Y \setminus P$ is open in $D(X \times Y)$. This proves~\eqref{eq:closed}.

On the other hand, write $0_X$ and $0_Y$ for the common origins in $X$ and $Y$, respectively. Then
\begin{equation}\label{eq:closure}
(0_X,0_Y) \in \cl P \setminus P \text{ in } D(X) \times D(Y).
\end{equation}
This is because every open neighborhood of $(0_X,0_Y)$ in $D(X) \times D(Y)$ contains some
\[
V := (\bigvee_s (-a_s,a_s)) \times (\bigvee_j (-b_j,b_j)).
\]
Now choose a sequence $t = (t_1,t_2,\ldots)$ with $t_j > \max\{j,1/b_j\}$ for all $j$, and choose an integer $k > 1/a_t$. Then $t_k > k > 1/a_t$ and hence $1/t_k < a_t$. Together with $1/t_k < b_k$, we have $(1/t_k,1/t_k) \in V \cap P$. This proves \eqref{eq:closure}.

The claim follows by combining~\eqref{eq:closed} and~\eqref{eq:closure}.
\end{ex}

\begin{rem}\label{rem:notcommute}
The same proof as above shows that $D(X \times X) \neq D(X) \times D(X)$, by replacing the positive integers $j$ to the constant sequences $(j,j,\ldots)$.
\end{rem}

We will systematically produce many examples of diffeological (vector) spaces $Z$ such that $D(Z \times Z) \neq D(Z) \times D(Z)$ in Theorem~\ref{thm:countable-vs-uncountable}.

\section{Free diffeological vector spaces vs free topological vector spaces}\label{s:comparison}

Recall from~\cite[proof of Proposition~3.5]{W} that for every diffeological space $X$, there exists a diffeological vector space $F(X)$ with the following universal property. Namely there exists a canonical smooth map $i_X:X \to F(X)$ such that for any diffeological vector space $V$ and any smooth map $f:X \to V$, there exists a unique smooth linear map $g:F(X) \to V$ with $f = g \circ i_X$. The diffeological vector space $F(X)$ is called the \textbf{free diffeological vector space} generated by the diffeological space $X$. It is constructed as $\oplus_{x \in X} \R$ as a vector space, equipped with the diffeology generated (as a diffeological space) by all maps
\[
\rho_n^X:(\R \times X)^n \to F(X)
\]
defined by
\[
(r_1,x_1,\ldots,r_n,x_n) \mapsto \sum_{i=1}^n r_i [x_i]
\]
for all $n \in \Z^{>0}$, where $[x_i]$ denotes a basis element in $F(X)$ corresponding to $x_i \in X$.

Given a topological space $Y$, one could ask for a topological vector space $G(Y)$ with a similar universal property. Namely, there exists a canonical continuous map $j_Y:Y \to G(Y)$ such that for any topological vector space $W$ and any continuous map $h:Y \to W$, there exists a unique continuous linear map $k:G(Y) \to W$ so that $h = k \circ j_Y$. When $Y$ is Tychonoff, it has been proven that $G(Y)$ exists; see e.g. \cite[Theorem~2.3(i)]{GM}. The proof there uses the abstract adjoint functor theorem. Note that $G(Y)$ as a vector space is also $\oplus_{y \in Y} \R$. \cite{Gp} gives a description of a basis for the topology of $G(Y)$.

We have the following simple observation connecting the two situations:

\begin{prop}\label{D(F(X)) = G(D(X))}
Let $X$ be a diffeological space such that $D(X)$ is Tychonoff. Then the topology of $G(D(X))$ is coarser than the topology of $D(F(X))$. If we further assume that $D(F(X))$ is a topological vector space\footnote{This is a key assumption that has to be checked independently in practice.}, then $D(F(X)) = G(D(X))$.
\end{prop}
In other words, under all of the above assumptions, $A \subseteq G(D(X))$ is open (resp. closed) if and only if for all $n \in \Z^{>0}$, $(\rho_n^X)^{-1}(A)$ is $D$-open (resp. $D$-closed) in the diffeological space $(\R \times X)^n$. The latter is equivalent to that for any $n \in \Z^{>0}$ and any plots $p_1:U_1 \to X, \cdots, p_n:U_n \to X$ of $X$, if we write $\rho_n^{p_1,\ldots,p_n}$ for the map $\R \times U_1 \times \cdots \times \R \times U_n \to F(X)$ defined by $(r_1,u_1,\ldots,r_n,u_n) \mapsto r_1 [p_1(u_1)] + \cdots + r_n [p_n(u_n)]$, then $(\rho_n^{p_1,\ldots,p_n})^{-1}(A)$ is open in the domain $\R \times U_1 \times \cdots \times \R \times U_n$. (According to the main results proved later on, we know that the assumptions in the above proposition hold (and so is the conclusion) when $X$ is an open subset of some Euclidean space (Theorem~\ref{thm:tvs}), and the second assumption and the conclusion both fail when $X$ is a wedge at $0$ of copies of $\R$ indexed over a continuum (Theorem~\ref{prop:surprise2}).) Compared with the concrete description of the topology, or the abstract universal property, this could be useful in practice.
\begin{proof}
We already know that the underlying vector spaces of $D(F(X))$ and $G(D(X))$ coincide.

Notice that the existence of the continuous map $Y \to G(Y)$ with codomain a topological vector space implies that for all $n \in \Z^{>0}$, the maps $\rho_n^Y:(\R \times Y)^n \to G(Y)$ defined as above must be continuous, when $Y$ is a Tychonoff (topological) space. In the current situation, since the identity map $D((\R \times X)^n) \to (\R \times D(X))^n$ of the underlying sets is continuous, this implies that the topology of $G(D(X))$ is coarser than the topology of $D(F(X))$, which is the first claim.

Use the universal property for $G(D(X))$ together with the further assumption in the proposition, we get a continuous linear map $G(D(X)) \to D(F(X))$. We conclude that the topology on $G(D(X))$ must be finer than the topology of $D(F(X))$. Combining with the first claim, we have $D(F(X)) = G(D(X))$, which is the second claim.
\end{proof}

As a comment, the map $j_Y:Y \to G(Y)$ is an inclusion of a closed (topological) subspace (\cite[Theorem~2.3(v)]{GM}) when $Y$ is Tychonoff, while $i_X:X \to F(X)$ can even fail to be an inclusion of a (diffeological) subspace (\cite[Example~3.7]{Wp}).

\section{$D$-open and $D$-closed subsets in $F(U)$}\label{sss:general}

Now we fix an arbitrary non-empty open subset $U$ of some Euclidean space\footnote{We reserve $U$ for an open subset of some Euclidean space throughout the paper.}, and we study the $D$-topology on the free diffeological vector space $F(U)$ generated by $U$.

Here are some (more or less) known results from the literatures:

We already know from~\cite[Corollary~6.4]{W} that every free diffeological vector space generated by a smooth manifold is projective, and from~\cite[Theorem~1.1]{CW19} that the $D$-topology on every projective diffeological vector space is Hausdorff. Therefore, every finite subset of $F(U)$ is $D$-closed.

Assume that a linear subspace $V$ of $F(U)$ is a smooth direct summand. (This is the case for every finite dimensional (resp. codimensional) linear subspace of $F(U)$ by~\cite[Theorems~1.1 and~4.2]{CW19}.) Write $\pi:F(U) \to W$ for the projection to the complement. As $D(W)$ is finer than the subspace topology of $D(F(U))$, the $D$-topology on $F(U)$ is Hausdorff, and $\pi^{-1}(0) = V$, we know that $V$ is $D$-closed in $F(U)$.

\begin{prop}
Let $A$ be a closed subset of $U$. Then the linear subspace of $F(U)$ spanned by the basis $A$ is $D$-closed in $F(U)$.
\end{prop}
\begin{proof}
This can be proved almost identically as that of~\cite[Lemma~2.4]{GM}.
\end{proof}

The following result is part of~\cite[Theorem~2.3(v)]{GM}, but the proof is adjusted to the diffeological case. The idea of the proof will be used frequently later, especially in Section~\ref{s:general}.

\begin{prop}
The image of the map $i_U:U \to F(U)$ sending $u$ to $[u]$ is $D$-closed in $F(U)$.
\end{prop}
\begin{proof}
Write $[U]$ for this image. For $\rho_n:(\R \times U)^n \to F(U)$, we need to show that $\rho_n^{-1}([U])$ is closed in $(\R \times U)^n$. Let $(r_1^j, u_1^j, \ldots, r_n^j, u_n^j)$ be a sequence in $\rho_n^{-1}([U])$ with respect to $j$ which converges to $(r_1^0, u_1^0, \ldots, r_n^0, u_n^0)$ in $(\R \times U)^n$. By definition and the assumption, for each $j$, $\sum_{i=1}^n r_i^j [u_i^j] = [u^j]$ for some $u^j\in U$. Without loss of generality, we may assume that $u_1^0 = \ldots = u_{k_1}^0 =: v_1$, $u_{k_1 + 1}^0 = \ldots = u_{k_2}^0 =: v_2$, $\cdots$, $u_{k_{l-1} + 1}^0 = \ldots = u_n^0 =: v_l$, and these $v_k$'s are mutually distinct. Hence, for $j$ sufficiently large, we know that the sets $\{u_1^j, \ldots, u_{k_1}^j\}$, $\cdots$, $\{u_{k_{l-1} + 1}^j, \ldots, u_n^j\}$ are mutually disjoint. Again without loss of generality, we may take a subsequence so that $\sum_{i = 1}^{k_1} r_i^j = 1$ and $\sum_{i = k_1 + 1}^{k_2} r_i^j = \ldots = \sum_{i = k_{l-1} + 1}^{n} r_i^j = 0$ for all $j$. This then implies that $\sum_{i = 1}^{k_1} r_i^0 = 1$ and $\sum_{i = k_1 + 1}^{k_2} r_i^0 = \ldots = \sum_{i = k_{l-1} + 1}^{n} r_i^0 = 0$. Hence $\rho_n(r_1^0, u_1^0, \ldots, r_n^0, u_n^0)
 =[v_1]\in [U]$, i.e., $(r_1^0, u_1^0, \ldots, r_n^0, u_n^0) \in \rho_n^{-1}([U])$. Therefore, $[U]$ is $D$-closed in $F(U)$.
\end{proof}

\begin{rem}\label{rem:closed}
The same argument shows that if $A \subseteq U$ is closed, then the image $[A]$ of the restriction to $A$ of the canonical map $i_U:U \to F(U)$ is also $D$-closed in $F(U)$. Here the closedness of $A$ is used to show that the image $\sum_i r_i^0 [u_i^0]$ of the limit $(r_1^0,u_1^0,\ldots,r_n^0,u_n^0)$ appeared in the domain in a similar proof to that of the above proposition is still in $A$.
\end{rem}

Now we take a systematic approach as follows:

For $m \in \N$, write $F_m(U)$ for the elements of $F(U)$ of the form $\sum_{i=1}^m r_i [u_i]$ with $r_i \neq 0$ and $u_i \neq u_j$ for all $i,j$ with $i \neq j$. Notice that $F_0(U) = \{0\}$. Write $F_{\leq n}(U) := \bigcup_{m \leq n} F_m(U)$. Therefore, we get a countable filtration of subspaces of $F(U)$:
$$
F_0(U) \subseteq F_{\leq 1}(U) \subseteq F_{\leq 2}(U) \subseteq \cdots \subseteq F(U).
$$

As a warm-up, we have:

\begin{prop}
Let $A$ be a subset of $F(U)$ such that $A \cap F_m(U)$ is finite for every $m \in \N$. Then $A$ is $D$-closed and $D$-discrete in $F(U)$.
\end{prop}
\begin{proof}
For any fixed $\rho_n:(\R \times U)^n \to F(U)$, notice that $\rho_n^{-1}(A) = \rho_n^{-1} \big{(} \bigcup_{m \leq n} (A \cap F_m(U) ) \big{)}$. By assumption, $\bigcup_{m \leq n} \big{(}A \cap F_m(U) \big{)}$ is a finite subset in $F(U)$, and hence $D$-closed. It implies that $\rho_n^{-1}(A)$ is closed in $(\R \times U)^n$. Therefore, $A$ is $D$-closed in $F(U)$.

For any $x \in A$, write $B$ for $A \setminus \{x\}$. Then the subset $B$ of $F(U)$ also has the property that $B \cap F_m(U)$ is finite for every $m \in \N$, since $B \cap F_m(U) \subseteq A \cap F_m(U)$. The above proved result shows that $B$ is $D$-closed in $F(U)$, which implies that $F(U) \setminus B$ is a $D$-open neighborhood of $x$ which contains no other element of $A$, and hence $A$ is $D$-discrete in $F(U)$.
\end{proof}

\begin{prop}\label{prop:closed}
Each subset $F_{\leq m}(U)$ is $D$-closed in $F(U)$.
\end{prop}
\begin{proof}
For $\rho_n:(\R \times U)^n \to F(U)$ with $(r_1,u_1,\ldots,r_n,u_n) \mapsto \sum_{i=1}^n r_i [u_i]$, we are left to show that $\rho_n^{-1}(F_{\leq m}(U))$ is closed in $(\R \times U)^n$.

For $n \leq m$, we know that $\rho_n^{-1}(F_{\leq m}(U)) = (\R \times U)^n$. So the result follows.

Now consider the case $n > m$. Take a sequence $(r_1^j,u_1^j,\ldots,r_n^j,u_n^j)$ with respect to $j$ in $\rho_n^{-1}(F_{\leq m}(U))$ which converges to $(r_1^0,u_1^0,\ldots,r_n^0,u_n^0)$ in $(\R \times U)^n$. Without loss of generality, we may assume that $u_1^0 = \ldots = u_{k_1}^0 =: v_1$, $u_{k_1 + 1}^0 = \ldots = u_{k_2}^0 =: v_2$, $\cdots$, $u_{k_{l-1} + 1}^0 = \ldots = u_n^0 =: v_l$, and these $v_k$'s are mutually distinct. It is enough to consider the case when $l > m$.  Notice that for $j$ sufficiently large, the sets $\{u_1^j, \ldots, u_{k_1}^j\}$, $\cdots$, $\{u_{k_{l-1} + 1}^j, \ldots, u_n^j\}$ are mutually disjoint. By definition of $F_{\leq m}(U)$, at most $m$-terms of $\sum_{i=1}^{k_1} r_i^j$, $\cdots$, $\sum_{i = k_{l-1} + 1}^n r_i^j$ can be non-zero. There are only finitely many such choices. By passing to a subsequence and rearranging, we may assume that the first $(l - m)$-terms of such sums are always $0$ for all $j$, which will lead to the same equality in the limit case. This then implies that $(r_1^0,u_1^0,\ldots,r_n^0,u_n^0) \in \rho_n^{-1}(F_{\leq m}(U))$, as we desired.
\end{proof}

\begin{rem}
Here is an alternative proof of the above proposition. Write $F_{\geq m+1}(U) := \bigcup_{i = m+1}^\infty F_i (U)$. Then to show that $F_{\leq m}(U)$ is $D$-closed in $F(U)$ is equivalent to that $F_{\geq m+1}(U)$ is $D$-open in $F(U)$. Let $\rho_n:(\R \times U)^n \to F(U)$. If $(r_1,u_1,\ldots,r_n,u_n) \in \rho_n^{-1}(F_{\geq m+1}(U))$, then there exists a partition of $\{1,2,\ldots,n\}$ according to the values of $u_i$'s, say $I_1 \cup \cdots \cup I_l = \{1,2,\ldots,n\}$, such that there exists $k$ with $m+1 \leq k \leq l$ satisfying $\sum_{j \in I_i} r_j \neq 0$ for $i \leq k$ and $\sum_{j \in I_i} r_j = 0$ for $i > k$. Therefore, we can take $\delta > 0$ so small that these balls $B_\delta(u_i)$'s do not intersect for taking one representative $i$ from each $I_1,\ldots,I_l$, and simultaneously $\sum_{j \in I_i} r_j' \neq 0$ for any $r_j' \in (r_j - \delta, r_j + \delta)$ with $i < k$. This implies the $D$-openness of $F_{\geq m+1}(U)$ in $F(U)$.
\end{rem}

Now we can inductively construct $D$-open neighborhoods of any point $x \in F(U)$ as follows. Assume that $x \in F_m(U)$. Take arbitrary open neighborhood $V_m^y$ in $(\R \times U)^m$ for any $y \in \rho_m^{-1}(x)$, and denote $V_m := \bigcup_{y \in \rho_m^{-1}(x)} \rho_m(V_m^y)$. Now inductively assume that $V_n$ has been constructed with $n \geq m$. Take arbitrary open neighborhood $V_{n+1}^{w,z}$ in $(\R \times U)^{n+1}$ for any $w \in V_n$ and any $z \in \rho_{n+1}^{-1}(w)$, and denote $V_{n+1} := \bigcup_{w \in V_n, z \in \rho_{n+1}^{-1}(w)} \rho_{n+1}(V_{n+1}^{w,z})$. Clearly, $V_n \subseteq V_{n+1}$. Write $V$ for $\bigcup_{n=m}^\infty V_n$. (Note that, as the proof in the above remark, one can see that as long as we take small enough open neighborhood in each step, $V \subseteq F_{\geq m}(U)$.) Then we have:

\begin{thm}
Let $V$ be constructed as above. Then $V$ is a $D$-open neighborhood of $x$ in $F(U)$, and every $D$-open neighborhood of $x$ contains some $V$ of this form.
\end{thm}
\begin{proof}
For $\rho_l:(\R \times U)^l \to F(U)$, we have to show that $\rho^{-1}_l(V)$ is open in $(\R \times U)^l$. Assume that $y := (r_1,u_1,\ldots,r_l,u_l) \in \rho^{-1}_l(V)$ and $\rho_l(y) \in V_n$. Write $s = \max\{n+1,l\}$, and denote $z := (r_1,u_1,\ldots,r_l,u_l,0,u,\ldots,0,u)$ for any fixed $u \in U$. Then $\rho_l(y) = \rho_s(z)$. By the construction of $V_s$ from $V_{s-1}$, we know that $\rho_s^{-1}(V_s)$ contains an open neighborhood of $z$ in $(\R \times U)^s$, which implies that $\rho_l^{-1}(V_s)$ (and hence $\rho_l^{-1}(V)$) contains an open neighborhood of $y$ in $(\R \times U)^l$, as we desired for the first claim.

The second claim follows from the fact that the $D$-topology on $F(U)$ is determined by all these $\rho_l$'s and the construction of $V$, together with the proof of the first claim.
\end{proof}

One might celebrate to have good understanding of the $D$-topology of $F(U)$ from the above theorem. Indeed, it is not so. The key difficulty is the subtle interrelationship between the possible sizes of the open neighborhoods of the preimage in $(\R \times U)^l$ of different expressions of the same element in $F(U)$. This will be overcome by the $k_\omega$-space theory in the following two sections.

\section{$k_\omega$-space theory}\label{s:komega}
In this section, we study the conservation of finite product (with the product topology) for a special class of $k$-spaces called $k_\omega$-spaces. Most results in this section are known to general topologists and they are scattered in the literatures; see~\cite{Michael-1966}, \cite{Franklin-Thomas}, \cite[Page 23]{Lin} and~\cite[Section~2.8]{Sakai-book-2013}. For the convenience of readers, we give detailed proofs. The authors would like to thank Professor Shou Lin for providing details for some parts of this section.

By a space here, we mean a topological space, since this section is purely topological. Recall that a $k$-space is a Hausdorff space such that a subset is closed if and only if its intersection with every compact subset is closed.

\begin{de}\
\begin{enumerate}
\item A Hausdorff space $X$ is called a \textbf{$k_\omega$-space} if there exists a countable cover $\{K_n:n\in\N\}$ (which is called a \textbf{$k_\omega$-structure} of $X$) of $X$ consisting of compact sets in $X$ such that a set $F$ is closed in $X$ if and only if $F\cap K_n$ is closed in $K_n$ for every $n\in\N$. If we further have $K_n \subseteq K_{n+1}$ for each $n$, we say that this $k_\omega$-structure is \textbf{towered}.
\item A continuous map $f:X\to Y$ from a space $X$ onto a space $Y$ is called a \textbf{compact-cover map} if for every compact set $F$ in $Y$, there exists a compact set $K$ in $X$ such that $f(K)=F$.
\end{enumerate}
\end{de}
It is straightforward to check that a Hausdorff space has a $k_\omega$-structure if and only if it has a towered $k_\omega$-structure.

\begin{ex}\label{ex:komega}
Every second-countable locally compact Hausdorff space is a $k_\omega$-space.
\end{ex}
\begin{proof}
It is a well-known result in topology that for such space $X$, there exists a family $\{K_i\}_{i \in \N}$ of compact subsets with $X = \bigcup_{i} K_i$ and $K_i \subseteq \Int K_{i+1}$ for each $i$; see for example~\cite[Proposition~A.60]{L}. Now let $F$ be a subset of $X$ such that $F \cap K_i$ is closed in $K_i$ for each $i$. Take $x \in \cl F $. Then there exists a sequence $(x_n)_{n \in \N}$ in $F$ converging to $x$, and $x \in K_i$ for some $i$. We may assume that $x_n \in \Int K_{i+1}$ for all $n$. Since $F \cap K_{i+1}$ is closed in $K_{i+1}$, we must have $x \in F$. Hence, $\cl F = F$, which means that $X$ is a $k_\omega$-space with $\{K_i\}_{i \in \N}$ as a $k_\omega$-structure.
\end{proof}

\begin{prop}\label{quotient-of-$k_omega$}
Let $q:X\to Y$ be a quotient map from a $k_\omega$-space $X$ onto a  Hausdorff space $Y$. Then $Y$ is also a $k_\omega$-space.
\end{prop}
\begin{proof} Let $\{K_n\}$ be a $k_\omega$-structure of $X$.
We check that  $\{q(K_n)\}$ is a $k_\omega$-structure of $Y$. For every set $F$ in $Y$, if $F\cap q(K_n)$ is closed in $q(K_n)$ for every $n$, then  $F\cap q(K_n)$ is closed in $Y$ and hence $q^{-1}(F\cap q(K_n))=q^{-1}(F)\cap q^{-1}( q(K_n))$ is closed in $X$.  Since
$q^{-1}(F)\cap K_n=q^{-1}(F)\cap q^{-1}( q(K_n))\cap K_n$, we have that $q^{-1}(F)\cap K_n$ is closed in $K_n$ for every $n$. It follows from $\{K_n\}$ being a $k_\omega$-structure of $X$ that $q^{-1}(F)$ is closed in $X$. Therefore, by $q$ being a quotient map, $F$ is closed in $Y$.
\end{proof}

\begin{lem}\label{containment}
Let $\{K_n:n \in \N\}$ be a towered $k_\omega$-structure on a $k_\omega$-space $X$. Then every compact subset $C$ of $X$ is contained in some $K_n$.
\end{lem}
\begin{proof}
If it were not the case, then $C \setminus K_n \neq \emptyset$ for every $n \in \N$. Take some $x_1 \in C \setminus K_1$. Since $\bigcup_{n \in \N} K_n = X$, it makes sense to write $n_2$ for the smallest $m \in \N$ such that $x_1 \in K_m$. Then take some $x_2 \in C \setminus K_{n_2}$, and let $n_3$ be the smallest $m \in \N$ such that $x_2 \in K_m$, and so on. In this way, we get an infinite subset $E : = \{x_i:i \in \N \}$ of $C$, which intersects every $K_n$ with only finitely many points, and hence closed in $X$ by definition of $k_\omega$-space. The same argument shows that every subset of $E$ is closed in $X$, and hence $E$ is discrete, which then conflicts with the compactness of $C$.
\end{proof}

\begin{prop}\label{quotient-is-compact-convex}
Let $q:X\to Y$ be a quotient map from a $k_\omega$-space $X$ onto a Hausdorff space $Y$. Then $q:X\to Y$ is a compact-cover map.
\end{prop}
\begin{proof} Let $\{K_n\}$ be a towered $k_\omega$-structure of $X$.
By Lemma \ref{quotient-of-$k_omega$} and its  proof, $Y$ is a $k_\omega$-space and  $\{q(K_n)\}$ is a $k_\omega$-structure of $Y$. For every compact set $C$ in $Y$, there exists some $n \in \N$ such that $C\subseteq q(K_n)$ by Lemma~\ref{containment}. Thus, $q^{-1}(C)\cap K_n$ is compact in $X$ and $q(q^{-1}(C)\cap K_n)=C$. It follows that  $q:X\to Y$ is a compact-cover map.
\end{proof}

\begin{prop}\label{product-of-compact-cover}
For $i=1,2$, let
$q_i:X_i\to Y_i$ be a compact-cover map from $X_i$ onto a Hausdorff space $Y_i$. Then the map $q_1\times q_2:X_1\times X_2\to Y_1\times
Y_2$ is a compact-cover  map.
\end{prop}
\begin{proof} Let $\pi_i:Y_1\times Y_2\to Y_i$  be the projections for $i=1,2$. Let $K$ be an arbitrary compact set in  $Y_1\times Y_2$. Then, for $i=1,2$,  $\pi_i(K)$ is compact in $Y_i$ and hence there exists a compact set $C_i$ in $X_i$ such that $q_i(C_i)=\pi_i(K)$.
Then $(q_1\times q_2)^{-1}(K)\cap (C_1\times C_2)$ is compact in $X_1\times X_2$ and $(q_1\times q_2)\big((q_1\times q_2)^{-1}(K)\cap (C_1\times C_2)\big)=K$.
Therefore, $q_1\times q_2:X_1\times X_2\to Y_1\times Y_2$ is a compact-cover  map.
\end{proof}

\begin{prop}\label{compact-cover-onto-$k$-is-quotient}
  Every compact-cover map onto a $k$-space is a quotient map.
\end{prop}
\begin{proof} Let $f:X\to Y$ be a compact-cover map from $X$ onto a $k$-space $Y$. For every compact set $C$, choose a compact set $K_C$ in  $X$ such that $f(K_C)=C$. Consider a  set $F$ in $Y$ satisfying $f^{-1}(F)$ is closed in $X$. Then, for every compact set $C$ in $Y$, $f^{-1}(F)\cap K_C$ is compact in $X$ and $f(f^{-1}(F)\cap K_C)=F\cap C$ is compact and hence is closed in $C$. It follows from $Y$ being a $k$-space that $F$ is closed in $Y$. Hence, $f:X\to Y$ is a quotient map.
\end{proof}

\begin{thm}\label{product-k-omega-is-k-omega}
 The finite product (with the product topology) of $k_\omega$-spaces is a $k_\omega$-space.
\end{thm}
\begin{proof}
It is enough to prove the binary product case. Let $X$ and $Y$ be two $k_\omega$-spaces. Then $X \times Y$ is Hausdorff. Now, let $\{K_n\}$ and $\{L_n\}$ be towered $k_\omega$-structures of $X$ and $Y$, respectively. We check below that $\{K_n\times L_n\}$ is a $k_\omega$-structure of $X\times Y$.

It is clear that $\bigcup_{n \in \N} (K_n \times L_n) = X \times Y$. Moreover, let $F$ be a set in $X\times Y$ such that $F\cap (K_n\times L_n)$ is closed for every $n$. We are left to show that $F$ is closed in $X \times Y$. For every $(x,y)\in (X\times Y)\setminus F$, choose $n$ such that $(x,y)\in K_n\times L_n$. We inductively define two sequences $\{U_i:i\geq n\}$ and $\{V_i:i\geq n\}$ of subsets of $X$ and $Y$, respectively, satisfying the following conditions for every $i\geq n$:

(i)~ $x\in U_i$, $y\in V_i$ and $U_i$ and $V_i$ are open in $K_i$ and $L_i$, respectively;

(ii)~ $(\cl U_i\times \cl V_i) \cap F=\emptyset$;

(iii)~$\cl U_i\subseteq K_i\cap U_{i+1}$ and  $\cl V_i\subseteq L_i\cap V_{i+1}$.\\
Since  $F\cap (K_n\times L_n)$ and $\{(x,y)\}$ are two disjoint closed subsets in the compact Hausdorff space $K_n \times L_n$, there exist open subsets $U_n$ and $V_n$ in $K_n$ and $L_n$, respectively, such that
$x\in U_n\subseteq \cl U_n\subseteq K_n$, $y\in V_n\subseteq \cl V_n\subseteq L_n$  and $(\cl U_n\times \cl V_n) \cap F=\emptyset$.
Note that  $\cl U_n\times \cl V_n$ is a compact set 
in the space $K_{n+1}\times L_{n+1}$ and $(\cl U_n\times \cl V_n) \cap F=\emptyset.$
It follows that there exist two open sets $U_{n+1}$ and $V_{n+1}$ in  $K_{n+1}$ and $L_{n+1}$, respectively, such that (i)-(iii) hold
for $n+1$.  Go on, and we obtain two such sequences.

Let $U=\bigcup_{i=n}^\infty U_i$ and $V=\bigcup_{i=n}^\infty V_i$.  Then $(x,y)\in U\times V$ and $(U\times V) \cap F=\emptyset.$
It remains to show that $U$ and $V$ are open in $X$ and $Y$. For example, we show the former. For every $i_0\geq n$, by (iii),
$$
U\cap K_{i_0}=\bigcup_{i=n}^\infty U_i\cap K_{i_0}=\bigcup_{i=i_0}^\infty U_i\cap K_{i_0}.
$$
Note that, for every $i\geq i_0$, $U_i\cap K_{i_0}$ is open in $K_{i_0}$. It follows from $\{K_n\}$ being a $k_\omega$-structure
of $X$ that $U$ is open in $X$.
\end{proof}
\begin{rem}
It is well-known that the product of two $k$-spaces is not necessarily a $k$-space, while Theorem~\ref{product-k-omega-is-k-omega} confirms the contrary for $k_\omega$-spaces. This will be extremely important for us later in the next sections.
\end{rem}

Recall that the category $\Top$ of topological spaces is bicomplete, i.e., limits and colimits exist. As an example,
let $\{X_n\}_n$ be a sequence of spaces such that $X_n$ is a subspace of $X_{n+1}$ for every $n$. Then we can define
a topology on the set $X=\bigcup_{n=1}^\infty X_n$ as follows: $A\subseteq X$ is open if and only if $A\cap X_n$ is open in
$X_n$ for every $n$, or equivalently, $A\subseteq X$ is closed if and only if $A\cap X_n$ is closed in
$X_n$ for every $n$. The space $X$ with this topology is called the \textbf{colimit} of  $\{X_n\}_n$ and denoted by $\underline{\rm colim}_n X_n$. Indeed, it is the colimit of the diagram
$$
 \xymatrix{
 X_1 \ \ar@{^{(}->}[r] & X_2 \ \ar@{^{(}->}[r] &\cdots \ \ar@{^{(}->}[r] & X_n \ \ar@{^{(}->}[r] &\cdots
 }$$
 in the category $\Top$.
By definition, every $k_\omega$-space is the colimit of its towered $k_\omega$-structure.

\begin{prop}\label{colim-of-normal-is-also}
Let $\{X_n:n \in \N\}$  be a sequence of $T_1$ normal spaces\footnote{There are different conventions for a normal space in the literatures, and we choose to express our statement in this way to avoid confusion.}, and let $X_n\subseteq X_{n+1}$ be closed for every $n$. Then $X=\underline{\rm colim}_n X_n$ is a $T_1$ normal space.
\end{prop}

\begin{proof}
The proof is essentially similar to that of Theorem~\ref{product-k-omega-is-k-omega}, so we will omit most of the details.

Trivially, $X$ is $T_1$. Now let $A$ and $B$ be two disjoint closed subsets of the space $X$. Using the fact that each $X_n$ is normal, inductively for $n \in \N$, we can construct open subsets $U_n$, $U_n'$, $V_n$ and $V_n'$ in $X_n$, and closed subsets $A_n,B_n$ in $X_n$, such that
\begin{enumerate}
\item $A_0 = A \cap X_0$ and $B_0 = B \cap X_0$;
\item $A_n \subseteq U_n \subseteq \cl U_n \subseteq U_n'$ and $B_n \subseteq V_n \subseteq \cl V_n \subseteq V_n'$;
\item $U_n' \cap V_n' = \emptyset$;
\item $A_{n+1} = (A \cap X_{n+1}) \bigcup \cl U_n$ and $B_{n+1} = (B \cap X_{n+1}) \bigcup \cl V_n$.
\end{enumerate}
Note that, for each $n$, $\cl_{X_n} S=\cl_X S$ for every $S\subseteq X_n$ since $X_n$ is closed in $X$. Define $U := \bigcup_{n \in \N} U_n$ and $V := \bigcup_{n \in \N} V_n$. It is straightforward to check that $U$ and $V$ are open neighborhoods of $A$ and $B$ respectively, such that $U \cap V = \emptyset$.
\end{proof}

\begin{rem}\label{colim-regulau-not-regulat}
Note that the colimit of a sequence of  $T_1$ regular spaces is not necessarily Hausdorff; see \cite[Section 2.9]{Sakai-book-2013}.
\end{rem}

\begin{prop}\label{k-omega-structure-of-k-omega}
 Let $\{K_n\}$  be  a  $k_\omega$-structure for a  $k_\omega$-space $X$, and let $\{X_n\}$  be  a  family of closed subsets in $X$. If $X_n\subseteq X_{n+1}$ and, for every $n$, there exists $m\in\N$ such that $K_n\subseteq X_m$, then $X=\underline{\rm colim}_n X_n$.
\end{prop}
\begin{proof} It is trivial.
\end{proof}

\begin{prop}\label{prop:normal}
Every $k_\omega$-space is normal.
\end{prop}
\begin{proof} It follows from Proposition \ref{colim-of-normal-is-also}.
\end{proof}
\begin{prop}\label{colim-of-$k_omega$-is-also}
 Let $\{X_n:n \in \N\}$  be a sequence of $k_\omega$-spaces, and let $X_n\subseteq X_{n+1}$ be closed for every $n$. Then $X=\underline{\rm colim}_n X_n$ is a $k_\omega$-space.
\end{prop}
\begin{proof} By Proposition \ref{colim-of-normal-is-also}, $X$ is Hausdorff. For every $n$, let $\{K_m^n:m\in\N\}$ be a  $k_\omega$-structure for $X_n$. We check that the family $\{K_m^n:n,m\in\N\}$ of compact sets in $X$ is a  $k_\omega$-structure for $X$. In fact, it is trivial to verify that it is a cover of $X$. Moreover, for some subset $F$ in $X$, assume that $F\cap K_m^n$ is closed in $K_m^n$ for each $n,m\in\N$.
Then, for every $n$, $F\cap X_n$ is closed since  $\{K_m^n:m\in\N\}$ is a  $k_\omega$-structure for $X_n$. Thus $F$ is closed in $X$ since $X$ is a colimit.
\end{proof}

\section{$k_\omega$-space theory applied to $F(U)$}\label{s:apply}

Now we use the $k_\omega$-space theory to show that $F(U)$ with the $D$-topology is a topological vector space (Theorem~\ref{thm:tvs}). Here $U$ is a non-empty open subset of some Euclidean space.

The use of $k_\omega$-space theory to $F(U)$ is based on the following observations:

Recall that in Section~\ref{sss:general}, we get a filtration of $F(U)$ by subspaces
$$
\{0\} = F_0(U) \subseteq F_{\leq 1}(U) \subseteq F_{\leq 2}(U) \subseteq \cdots \subseteq F(U),
$$
where $F_{\leq n}(U)$ denotes the set of elements in $F(U)$ which can be expressed as a sum of at most $n$-terms under the base $U$. We can equip each $F_{\leq n}(U)$ with the subset diffeology of $F(U)$.

On the other hand, for each $n$ we have a surjective map
\[
\rho_n:(\R \times U)^n \to F_{\leq n}(U)
\]
given by $(r_1,u_1,\ldots,r_n,u_n) \mapsto \sum_{i=1}^n r_i [u_i]$. So we can equip the set $F_{\leq n}(U)$ with the quotient diffeology from $\rho_n$. This diffeology in general differs from the subset diffeology of $F(U)$ due to cancellation, so let us write this new diffeological space as $F'_{\leq n}(U)$. It follows immediately from definition of colimit in the category $\Diff$ of diffeological spaces that
\begin{equation}\label{eq:diffeomorphism}
F(U) \cong \underline{\rm colim}_n F_{\leq n}(U) \cong \underline{\rm colim}_n F'_{\leq n}(U)
\end{equation}
as diffeological spaces, where the inclusion $F_{\leq n}'(U) \subseteq F_{\leq n+1}'(U)$ of the underlying sets is realized as a smooth map $F_{\leq n}'(U) \to F_{\leq n+1}'(U)$ induced by $(\R \times U)^n \to (\R \times U)^{n+1}$ with $(r_1,u_1,\ldots,r_n,u_n) \mapsto (0,u,r_1,u_1,\ldots,r_n,u_n)$ for any fixed $u \in U$. We do not know in priori the relationship between $D(F_{\leq n}'(U))$ and $D(F_{\leq n+1}'(U))$.

Although the diffeological structures on the set $F_{\leq n}(U)$ discussed above can differ, their induced $D$-topologies coincide:

\begin{lem}\label{subspace}
For every $n$, $D(F_{\leq n}(U)) = D(F'_{\leq n}(U))$, and $D(F'_{\leq n}(U))$ is a closed subspace of $D(F'_{\leq n+1}(U))$.
\end{lem}
\begin{proof}
Since the $D$-topology functor is a left adjoint, the diffeomorphism~\eqref{eq:diffeomorphism} above implies that
$$
D(F(U)) \cong \underline{\rm colim}_n D(F'_{\leq n}(U)).
$$
By Proposition~\ref{prop:closed}, to prove the first statement, it is enough to prove the second one, i.e., $D(F'_{\leq n}(U))$ is a closed subspace of $D(F'_{\leq n+1}(U))$.

Take any fixed $u_0 \in U$, the smooth map $(\R \times U)^n \to (\R \times U)^{n+1}$ given by $(r_1,u_1,\ldots,r_n,u_n) \mapsto (0,u_0,r_1,u_1,\ldots,r_n,u_n)$ induces a smooth inclusion $F'_{\leq n}(U) \to F'_{\leq n+1}(U)$, which implies that the $D$-topology on $F'_{\leq n}(U)$ contains the sub-topology from the $D$-topology of $F'_{\leq n+1}(U)$. For the converse inclusion, let $E$ be a closed subset in $D(F'_{\leq n}(U))$. So $\rho_n^{-1}(E)$ is closed in $(\R\times U)^n$, and we are left to show that $\rho_{n+1}^{-1}(E)$ is closed in $(\R \times U)^{n+1}$. Let $(r_1^m,u_1^m,\ldots,r_{n+1}^m,u_{n+1}^m)$ be a sequence (with respect to $m$) in $\rho_{n+1}^{-1}(E)$ converging to $x_{n+1} := (r_1,u_1,\ldots,r_{n+1},u_{n+1})$. Then $E \subseteq F'_n(U)$ implies that for each fixed $m$, either some $r_i^m = 0$ or some $u_i^m = u_j^m$. There are only finitely many such choices. So by passing to a subsequence, we may assume that one of such cases occurs for all $m$, which then implies that the same equality holds in the limit and hence $x_{n+1} \in \rho_{n+1}^{-1}(E)$. It follows that $E$ is closed in $D(F'_{\leq n+1}(U))$.  Therefore, $D(F'_{\leq n}(U))$ is a closed subspace of $D(F'_{\leq n+1}(U))$.
\end{proof}

For the following two results, let $U,V$ be open subsets of some Euclidean spaces.

\begin{lem} \label{finite-continuity} For every $n$ and $m$, $D(F_{\leq n}(U))$ (resp. $D(F_{\leq m}(V))$) is a $k_\omega$-space and
$$
D(F'_{\leq n}(U)\times F'_{\leq m}(V)) = D(F_{\leq n}(U)) \times D(F_{\leq m}(V)).
$$
\end{lem}
\begin{proof}
By Lemma \ref{subspace} and the left adjointness of the $D$-topology functor, $\rho_n^U:(\R\times U)^n\to D(F_{\leq n}(U)) = D(F'_{\leq n}(U))$ and $\rho_m^V:(\R\times V)^m\to D(F_{\leq m}(V)) = D(F'_{\leq m}(V))$ are both quotient maps. Note that $(\R\times U)^n$ is a $k_\omega$-space by Example~\ref{ex:komega}. Since $D(F_{\leq n}(U))$ is Hausdorff, using Proposition~\ref{quotient-is-compact-convex}, $\rho_n^U:(\R\times U)^n\to D(F_{\leq n}(U))$ is a compact-cover map, and so is $\rho_m^V$. It follows from Proposition~\ref{product-of-compact-cover} that $\rho_n^U\times\rho_m^V:(\R\times U)^n\times (\R\times V)^m\to D(F_{\leq n}(U))\times D(F_{\leq m}(V))$ is a compact-cover.
By Proposition~\ref{quotient-of-$k_omega$} $D(F_{\leq n}(U))$ is a $k_\omega$-space, and so is $D(F_{\leq m}(V))$. Moreover, using Theorem~\ref{product-k-omega-is-k-omega}, we have that $D(F_{\leq n}(U))\times D(F_{\leq m}(V))$ is a $k_\omega$-space and hence a $k$-space. By Proposition~\ref{compact-cover-onto-$k$-is-quotient}, $\rho_n^U\times\rho_m^V:(\R\times U)^n\times (\R\times V)^m\to D(F_{\leq n}(U)) \times D(F_{\leq m}(V))$ is a quotient map. On the other hand, both $\rho_n^U:(\R \times U)^n \to F'_{\leq n}(U)$ and $\rho_m^V$ are subductions (i.e., equivalent to quotient maps in the category $\Diff$ of diffeological spaces), and hence so is $\rho_n^U \times \rho_m^V:(\R \times U)^n \times (\R \times V)^m \to F'_{\leq n}(U) \times F'_{\leq m}(V)$. Using the left adjointness of the $D$-topology functor, we get another quotient map $(\R \times U)^n \times (\R \times V)^m \to D(F'_{\leq n}(U) \times F'_{\leq m}(V))$. In summary, we get two quotient maps for the same map from the same topological space to the same set, and hence the two topologies on the codomain must coincide, which is the result.
\end{proof}

\begin{prop}\label{continuity} $D(F(U))$ and $D(F(V))$ are $k_\omega$-spaces, and
$$
D(F(U) \times F(V)) = D(F(U)) \times D(F(V)).
$$
\end{prop}
The result actually holds for finite product.
\begin{proof} By Proposition~\ref{colim-of-$k_omega$-is-also}, Lemma~\ref{finite-continuity} and the left adjointness of the $D$-topology functor,  both $D(F(U))=\underline{\rm colim}_n D(F_{\leq n}(U))$ and $D(F(V))$
are $k_\omega$-spaces. It follows from Theorem~\ref{product-k-omega-is-k-omega} that $D(F(U))\times D(F(V))$ is also a $k_\omega$-space. Moreover,  every compact set in $D(F(U))\times D(F(V))$ is contained in some $F_{\leq n}(U)\times F_{\leq n}(V)$ by Lemma~\ref{containment}. Therefore, by Lemma~\ref{k-omega-structure-of-k-omega},
\[
D(F(U))\times D(F(V))=\underline{\rm colim}_n \big( D(F_{\leq n}(U))\times D(F_{\leq n}(V)) \big).
\]
Using Lemma \ref{finite-continuity} and the left adjointness of the $D$-topology functor, we have
\[
D(F(U)) \times D(F(V)) = \underline{\rm colim}_n D(F'_{\leq n}(U) \times F'_{\leq n}(V)) = D \big( \underline{\rm colim}_n (F'_{\leq n}(U) \times F'_{\leq n}(V)) \big).
\]
Notice that the subsystem $\{F'_{\leq n}(U) \times F'_{\leq n}(V)\}_n$ is final in the system $\{F'_{\leq n}(U) \times F'_{\leq m}(V)\}_{n,m}$, the result then follows from the cartesian closedness of the category $\Diff$:
\begin{equation*}
\begin{split}
\underline{\rm colim}_n (F'_{\leq n}(U) \times F'_{\leq n}(V)) & = \underline{\rm colim}_{n,m} (F'_{\leq n}(U) \times F'_{\leq m}(V)) \\ & = \underline{\rm colim}_n (F'_{\leq n}(U) \times \underline{\rm colim}_m F'_{\leq m}(V)) \\ & = \underline{\rm colim}_n (F'_{\leq n}(U) \times F(V)) \\ & = (\underline{\rm colim}_n F'_{\leq n}(U)) \times F(V) \\ & = F(U) \times F(V),
\end{split}
\end{equation*}
as desired.
\end{proof}

As a digression, the same proof as above gives the following result which will be used in the next section.

\begin{prop}\label{prop:general-result}
Assume that $X = \underline{\rm colim}_n X_n$ for a countable sequence of injective smooth maps
\[
\ldots \to X_n \to X_{n+1} \to \ldots
\]
of diffeological spaces such that each image of $X_n$ is $D$-closed in $X_{n+1}$, each $D(X_n)$ is a $k_\omega$-space, and $D(X_n \times X_n) = D(X_n) \times D(X_n)$. Then $D(X)$ is also a $k_\omega$-space with $D(X \times X) = D(X) \times D(X)$.
\end{prop}

Let us go back to the theme of this section.

\begin{cor}\label{cor:addition}
The addition $+:D(F(U))\times D(F(U))\to D(F(U))$ is continuous.
\end{cor}
\begin{proof}
This follows from Proposition~\ref{continuity} by taking $U=V$ and the smoothness of the addition $F(U) \times F(U) \to F(U)$.
\end{proof}

\begin{thm}\label{thm:tvs}
For every open subset $U$ of some Euclidean space, the $D$-topology on $F(U)$ makes it a normal topological vector space.
\end{thm}
\begin{proof}
This follows from Corollary~\ref{cor:addition} together with Propositions~\ref{prop:normal} and~\ref{continuity}.
\end{proof}

\section{The $D$-topology on projective diffeological vector spaces}\label{s:projective}

Thanks to Theorem~\ref{thm:tvs} and the structural theorem of projective diffeological vector space (\cite[Corollary~6.15]{W}), we show in this section that every countably-generated projective diffeological vector space with the $D$-topology is a topological vector space.

We first recall from~\cite{W}:

\begin{de}\cite[Definition~6.1]{W}
A diffeological vector space $P$ is called \textbf{projective}, if for every linear subduction $V \to W$ and every smooth linear map $P \to W$, there exists a smooth linear map $P \to V$ making the following triangle commutative:
\[
\xymatrix{& P \ar[d] \ar@{.>}[dl] \\ V \ar[r] & W.}
\]
\end{de}

Since not every linear subduction has a smooth linear section, projective diffeological vector spaces provide the basis for doing homological algebra on the category of diffeological vector spaces. Here is the structural theorem for projective diffeological vector spaces:

\begin{thm}\label{thm:structural}\cite[Corollary~6.15]{W}
Every projective diffeological vector space is a smooth direct summand of a coproduct of free diffeological vector spaces generated by open subsets of Euclidean spaces.
\end{thm}

Projective diffeological vector spaces behave nicely. For example, \cite{CW19} showed that they are always Hausdorff under the $D$-topology, etc.

To reach our main result in this section, we will first prove the following results:

\begin{prop}\label{prop:directsum}
Let $\{U_i\}_{i \in I}$ be a set of open subsets of Euclidean spaces. Then
\[
D(\bigoplus_{i \in I} F(U_i)) = \bigoplus_{i \in I} D(F(U_i)),
\]
where $\bigoplus_{i \in I} D(F(U_i))$ is the direct sum of the family $\{ D(F(U_i) : i \in I)\}$
of topological vector spaces in the category $\mathfrak{T}\mathrm{VS}$ of topological vector spaces.
\end{prop}
\begin{proof}
When the index set $I$ is finite, the result follows from generalization of Proposition~\ref{continuity} to the finite product case. So for the rest of the proof, we may assume $I$ is infinite. Note that
\[
\bigoplus_{i \in I} F(U_i) = \underline{\rm colim}_J \bigoplus_{j \in J} F(U_j)
\]
as diffeological spaces, where $J$ runs over all finite subsets of $I$ (\cite[Proposition~3.2]{W}). Using the left adjointness of the $D$-topology functor together with the generalization of Proposition~\ref{continuity} to the finite product case, we have
\[
D(\bigoplus_{i \in I} F(U_i)) = \underline{\rm colim}_J D(\bigoplus_{j \in J} F(U_j)) = \underline{\rm colim}_J (\bigoplus_{j \in J} D(F(U_j))) = \bigoplus_{i \in I} D(F(U_i)),
\]
where the last equality is a result of topological vector spaces which can be deduced directly.
\end{proof}

\begin{prop}\label{prop:countable}
Let $V$ be a countable direct sum of $F(U_i)$'s with each $U_i$ open subset of some Euclidean space. Then $V$ is a topological vector space and a $k_\omega$-space under the $D$-topology with $D(V \times V) = D(V) \times D(V)$.
\end{prop}
From Corollary~\ref{prop:surprise1}, the countability condition of this proposition can \emph{not} be removed.
\begin{proof}
It is enough to prove the infinite countable case. Write $V = \bigoplus_{i=1}^\infty F(U_i)$. Then $V = \underline{\rm colim}_{j=1}^\infty V_j$ where $V_j = \bigoplus_{i=1}^j F(U_i) = \prod_{i=1}^j F(U_i)$ is a Hausdorff topological vector space under the $D$-topology with $D(V_j \times V_j) = D(V_j) \times D(V_j)$. Since $V_j$ is a direct summand of $V_{j+1}$, it is a $D$-closed subspace. Moreover, each $V_j$ is a $k_\omega$-space. Proposition~\ref{prop:general-result} shows that $D(V)$ is a $k_\omega$-space and $D(V \times V) = D(V) \times D(V)$, and hence the result.
\end{proof}

\begin{prop}\label{prop:directsummand}
Let $V$ be a diffeological vector space, and let $W$ be a smooth direct summand of $V$. Assume that $D(V \times V) = D(V) \times D(V)$. Then $D(W)$ is also a topological vector space with $D(W \times W) = D(W) \times D(W)$.
\end{prop}
\begin{proof}
The condition in the proposition implies that $D(V)$ is a topological vector space. We have a commutative diagram
\[
\xymatrix{D(W \times W) \ar[r] \ar[d] & D(V \times V) \ar[r] \ar[d] & D(W \times W) \ar[d] \\ D(W) \times D(W) \ar[r] & D(V) \times D(V) \ar[r] & D(W) \times D(W),}
\]
in the category of topological spaces, where each verticle map is induced by the two projections, the left horizontal maps are induced by the inclusion $W \to V$, the right horizontal maps are induced by a fixed retract $V \to W$, and the compositions of both horizontal maps are the identity maps. Therefore, $D(W \times W) = D(W) \times D(W)$, and $D(W)$ is a topological vector space.
\end{proof}

A similar retraction argument as above shows:

\begin{prop}\label{prop:retract}
Let $V$ be a diffeological vector space, and let $W$ be a smooth direct summand of $V$. If $D(V)$ is a topological vector space, then so is $D(W)$.
\end{prop}

We say that a projective diffeological vector space is \textbf{countably-generated} if it is a smooth direct summand of a countable direct sum of $F(U)$'s with each $U$ open subset of some Euclidean space.

As a direct conclusion, we have:

\begin{thm}\label{thm:proj}
For every countably-generated projective diffeological vector space $P$, $D(P)$ is a $k_\omega$-space (hence normal) and a topological vector space with $D(P \times P) = D(P) \times D(P)$.
\end{thm}
We will show in next section that this theorem could fail if the projective diffeological vector space is not countably-generated. More precisely and surprisingly, in this case, the projective diffeological vector space can fail to be a topological vector space under the $D$-topology.
\begin{proof}
This follows immediately from definition of countably-generated projective diffeological vector space, and Propositions~\ref{prop:countable} and~\ref{prop:directsummand}.
\end{proof}

\begin{cor}\label{cor:proj}
For every smooth manifold $M$, $D(F(M))$ is a topological vector space with $D(F(M) \times F(M)) = D(F(M)) \times D(F(M))$.
\end{cor}
More generally, if a diffeological space $X$ is countably-generated (i.e., it has a countable generating set of plots with the union of the images being $X$) and $F(X)$ is projective, then $F(X)$ is a topological vector space under the $D$-topology. For example, this applies to $X$ being the cross with the gluing diffeology\footnote{This also applies to $X$ being countably many $\R$'s glued at the origins. Compare this with Theorem~\ref{prop:surprise2}.}, or the orbifolds $\R^n/\pm 1$ or $\R^2/C_n$; see~\cite{W}. 
\begin{proof}
From Theorem~\ref{thm:proj}, it is enough to show that $F(M)$ is countably-generated. From our convention, every smooth manifold is assumed to be Hausdorff and second-countable, so one has a countable atlas for $M$, which implies that $F(M)$ is countably-generated.
\end{proof}

\section{Examples and further conclusions}\label{s:general}

In this section, we give two positive results and a negative result. More precisely, we show in Example~\ref{ex 9.1} and Proposition~\ref{prop:irrational} that two non-projective diffeological vector spaces are also topological vector spaces under the $D$-topology, and we show that in general the $D$-topology does not make a diffeological vector space (which is even both free and projective) into a topological vector space (Theorem~\ref{prop:surprise2}), based on the results in Section~\ref{s:notcommute}. From the latter, more conclusions are derived.

\begin{ex}\label{ex 9.1}
From~\cite[Example~4.3]{W}, we know that the countable product $\prod_\omega \R$ is not a projective diffeological vector space. By Claim~2 there, we know that $D(\prod_\omega \R) = \prod_\omega D(\R)$, and hence $D(\prod_\omega \R \times \prod_\omega \R) = D(\prod_\omega \R) \times D(\prod_\omega \R)$. So the $D$-topology makes $\prod_\omega \R$ into a topological vector space. 
\end{ex}

Irrational tori are key examples in diffeology: on every such space, the $D$-topology is indiscrete (i.e., the only open subsets are the empty set and the whole space), but one can still do some geometry on such spaces using diffeology; see~\cite{DI}. An irrational torus is defined to be the quotient group (equipped with the quotient diffeology) of the usual $2$-dimensional torus modulo an irrational line, or equivalently, it is the quotient group $\R/(\Z + \alpha \Z) =: T_\alpha$ of $\R$ (equipped with the quotient diffeology) with $\alpha \in \R \setminus \Q$. It is known from~\cite[Example~6.9(2)]{W} or~\cite[Corollary~3.17]{CW19} that the free diffeological vector space $F(T_\alpha)$ generated by the irrational torus $T_\alpha$ is not projective.

Note that unlike $D(T_\alpha)$, $D(F(T_\alpha))$ is not indiscrete. In fact, this is the case for every free diffeological vector space $F(X)$ generated by a diffeological space $X$, since the constant map $X \to \R$ with image $1 \in \R$ induces a smooth linear map $F(X) \to \R$, which produces many non-trivial $D$-open subsets for $F(X)$.

\begin{prop}\label{prop:irrational}
$F(T_\alpha)$ is a topological vector space under the $D$-topology.
\end{prop}
\begin{proof}
We write $\rho:F(T_\alpha) \to \R$ for the smooth linear map induced by the constant map $T_\alpha \to \R$ with $x \mapsto 1$:

Claim~1: $(T_\alpha)^n$ has indiscrete $D$-topology for all $n \in \N$.

This follows immediately from the fact that the preimage of any point in $(T_\alpha)^n$ for the quotient map $\R^n \to (T_\alpha)^n$ is dense in $\R^n$.

Claim~2: For every $D$-open subset $A$ of $F(T_\alpha)$, we have that $\rho^{-1}(\rho(A))=A$.

To prove this, we only need to show that if $\sum_{i=1}^n r_i [a_i] \in A $, then $\rho^{-1}(\sum_{i} r_i) \subseteq A$. Note that Claim~1 and the continuity of the smooth map $\rho_{r_1,\ldots,r_n}:(T_\alpha)^n \to F(T_\alpha)$ defined by $(x_1,\ldots,x_n) \mapsto \sum_{i=1}^n r_i[x_i]$ under the $D$-topology imply that its whole image is in $A$. In particular, for any fixed $x \in T_\alpha$, $(\sum_{i=1}^n r_i) [x] \in A$. So for any $(s_1,\ldots,s_m) \in \R^m$ with $\sum_{j=1}^m s_j = \sum_{i=1}^n r_i$, the image of $\rho_{s_1,\ldots,s_m}:(T_\alpha)^m \to F(T_\alpha)$ is in $A$. In other words, $\rho^{-1}(\sum_{i} r_i) \subseteq A$.

Claim~3: The $D$-topology of $F(T_\alpha)$ is determined by the map $\rho:F(T_\alpha) \to \R$. In other words, $B \subseteq F(T_\alpha)$ is $D$-open if and only if there exists an open subset $U$ of $\R$ such that $B = \rho^{-1}(U)$.

By Claim~2, we only need to verify that  $\rho(A)$ is open in $\R$ for every $D$-open subset $A$ in $F(T_\alpha)$.
For a fixed $x\in T_\alpha$, consider the $D$-continuous map $k:\R\to F(T_\alpha)$
given $r\mapsto r[x]$. It follows from Claim ~2 that $\rho(A)=k^{-1}(A)$. Hence, $\rho(A)$ is open in $\R$, as our required.

Claim~4: $D(F(T_\alpha) \times F(T_\alpha)) = D(F(T_\alpha)) \times D(F(T_\alpha))$.

From Claim~3, we know that the topology of the right-hand side is determined by the smooth map $\rho \times \rho: F(T_\alpha) \times F(T_\alpha) \to \R \times \R$, and we are left to show that so is the left-hand side. The verification goes almost identically as that of Claims~ 2 and 3.

The proposition then follows from Claim~4.
\end{proof}

\begin{rem}\
\begin{enumerate}
\item In fact, we can show a more general result that if $X$ is a diffeological space such that $D(X)$ is indiscrete, then $F(X)$ is a topological vector space under the $D$-topology. The proof goes almost identically as that of the previous proposition, except that of Claim~1, which can be proved as follows. Notice that for every non-trivial subset $B$ of $X^n$, (i.e., $B$ is neither $X^n$ nor the empty set), there exist $i$ and $(x_1,\ldots,x_{i-1},x_{i+1},\ldots,x_n) \in X^{n-1}$ such that the projection of $B \cap (\{x_1\} \times \cdots \{x_{i-1}\} \times X \times \{x_{i+1}\} \times \cdots \times \{x_n\})$ to the $i^{th}$ coordinate is a non-trivial subset of $X$. Then the smooth map $X \to X^n$ defined by $x \mapsto (x_1,\ldots,x_{i-1},x,x_{i+1},\ldots,x_n)$ implies that $B$ is not $D$-open in $X^n$. In other words, $D(X^n)$ is indiscrete.

\item Moreover, a similar Claim~3 as that in the proof of Proposition~\ref{prop:irrational} shows that $D(X)$ is a topological subspace of $D(F(X))$ which is not closed, via the canonical map $i_X:X \to F(X)$. Compare this with Remark~\ref{rem:Tychonoff}.

\item We can find a diffeological vector space which is not a free diffeological vector space such that it is a topological vector space under the $D$-topology. For example, take $V$ to be the vector space $\R$ with the indiscrete diffeology. Then it is clear that $V$ is not free, and $D(V)$ is indiscrete\footnote{Hence, $V$ is not projective.}. Then the addition $D(V) \times D(V) \to D(V)$ is automatically continuous, making $D(V)$ a topological vector space.
\end{enumerate}
\end{rem}

In the rest of this section, we will focus on the free diffeological vector space (which is also projective) generated by one of the diffeological spaces introduced in Example~\ref{ex:notcommute}, and show that it is not a topological vector space under the $D$-topology.

As the same as in Example~\ref{ex:notcommute}, let $X$ be the wedge at $0$ of copies of $\R$ indexed over all infinite sequences $s = (s_1,s_2,\ldots)$ of positive integers. It is straightforward to check that $D(X)$ is Hausdorff. Write $V$ for the free diffeological vector space generated by $X$. A similar argument as~\cite[Example~6.7]{W} shows that $V$ is a projective diffeological vector space.

\begin{rem}
One can also use homotopy theory developed in~\cite{CW14} to prove the projectivity of $V$ by the cofibrancy of $X$ as in~\cite[Remark~6.8]{W}. There is a typo in the statement of~\cite[Remark~6.8]{W}, and the correct one is as follows. Let $A$ and $B$ be diffeological spaces. If either $A$ is cofibrant, or $F(A)$ is a projective diffeological vector space and there is a cofibration $A \to B$ in the sense of~\cite{CW14}, then $F(B)$ is also a projective diffeological vector space.
\end{rem}

By cartesian closedness of $\Diff$, the diffeology on $V$ is generated by
\begin{equation}\label{eq:generate}
\rho_{s^1,\ldots,s^n}:\R \times \R_{s^1} \times \cdots \times \R \times \R_{s^n} \to V
\end{equation}
with $(r_1,x_{1,s^1},\ldots,r_n,x_{n,s^n}) \mapsto \sum_i r_i [x_{i,s^i}]$ for all $n \in \N$, where each $s^i$ is an infinite sequence of positive integers. Note that $[x_{i,s^i}] = [x_{j,s^j}]$ if and only if either $x_{i,s^i} = x_{j,s^j} \neq 0$ and $s^i = s^j$ or $x_{i,s^i} = x_{j,s^j} = 0$.

\begin{lem}\label{lem:Tychonoff}
$X$ is a \textbf{smooth Tychonoff (diffeological) space}, in the sense that for every point and every disjoint closed subset of $D(X)$, there exists a smooth map $X \to \R$ which sends the point to $0$ and the closed subset to $1$.
\end{lem}
\begin{proof}
Write $x$ for the given point. This lemma can be proved by two cases:

If $x$ is the common zero, then the complement of the closed subset contains an open interval $(-a_s,a_s)$ on each $\R_s$ inside $X$. One can define a smooth map $\R_s \to \R$ sending $0 \in \R_s$ to $0$ and $(-\infty,-a_s] \cup [a_s,\infty)$ to $1$ for each $s$. These together produce a smooth map $X \to \R$ with the required property.

If $x$ is not the common zero, then the complement of the closed subset contains an open interval in the copy $\R_s$ containing $x$ so that a neighborhood of $0 \in \R_s$ is apart from the interval. Then any smooth map $\R_s \to \R$ which sends $x$ to $0$ and the outside of the open interval to $1$ together with the rest part of $X$ to $1$, gives a smooth map $X \to \R$ with the required property.
\end{proof}

\begin{lem}
$D(X)$ is a topological subspace of $D(V)$ induced by the canonical map $i_X:X \to V$.
\end{lem}
\begin{proof}
It is enough to show that for every proper closed subset $A$ of $D(X)$, $[A]$ is closed in the topological subspace $[X]$ of $D(V)$, where $[A]$ (resp. $[X]$) denotes the image of the restriction to $A$ (resp. $X$) of the canonical map $i_X:X \to V$. Take any $x \in X \setminus A$. By Lemma~\ref{lem:Tychonoff}, we get a smooth map $\tau_x:X \to \R$ such that $x \mapsto 0$ and $A \mapsto 1$. By the universal property of free diffeological vector space, we get a smooth linear map $\tilde{\tau}_x:V \to \R$ such that $U_x := \tilde{\tau}_x^{-1}(\R \setminus \{1\})$ is open in $D(V)$ containing $[x] \in V$, and $[A]\cap U_x=\emptyset$. Write $U := \cup_{x \in X \setminus A} U_x$. Then $U$ is open in $D(V)$, and $(V \setminus U) \cap [X] = [A]$, which concludes the claim.
\end{proof}

\begin{rem}\label{rem:Tychonoff}
In fact, one can prove a stronger result than the previous lemma, but it is not needed later: Let $Y$ be a smooth Tychonoff space such that $D(Y)$ is $T_1$. Then $D(Y)$ is a closed subspace of $D(F(Y))$ via the canonocal map $i_Y: Y \to F(Y)$.

Notice that the proof of the previous lemma only uses definition of smooth Tychonoff space, so the conclusion holds for $Y$ and $F(Y)$. Hence we are left to show that $[Y]$ is closed in $D(F(Y))$. It suffices to show that for every $z \in F(Y) \setminus [Y]$, there exists a $D$-open neighborhood $W$ of $z$ in $F(Y)$ such that $W \cap [Y] = \emptyset$. Notice that if $z = \sum_{i=1}^n r_i [y_i] \in F_n(Y)$, then either $\sum_{i=1}^n r_i \neq 1$ or $\sum_{i=1}^n r_i = 1$ and $n \geq 2$. The collection of all $z$'s in the first case form a $D$-open subset of $F(Y)$ which is disjoint with $[Y]$, by definition of smooth Tychonoff. For $z$ in the second case, there exists $i$ such that $r_i \neq 1$. Since $D(Y)$ is $T_1$, $\{y_i\}$ and the closed subset $A := \{y_1,\ldots,y_n\} \setminus \{y_i\}$ of $D(Y)$ are disjoint. So there exists a smooth map $Y \to \R$ with $y_i \mapsto 1$ and $A \mapsto 0$, which induces a smooth linear map $F(Y) \to \R$ with $z \mapsto \R \setminus \{1\}$. Therefore, we proved the result.
\end{rem}

Therefore, we have:

\begin{cor}\label{cor:subspace}
$D(X) \times D(X)$ is a topological subspace of $D(V) \times D(V)$ induced by the canonical map $i_X:X \to V$.
\end{cor}

\begin{thm}\label{prop:surprise2}
For the diffeological vector space $V$ defined above (which is both free and projective), $D(V)$ is not a topological vector space.
\end{thm}
\begin{proof}
We will show that $f:D(V) \times D(V) \to D(V)$ induced by the addition $+:V \times V \to V$ is not continuous.

Note that $g:D(V \times V) \to D(V)$ induced by the addition $+:V \times V \to V$ is continuous, and all of the maps $f,g,+$ are the same as the underlying set maps. Take the set $P$ as in Remark~\ref{rem:notcommute}. More precisely, write $p_{s,j} = (1/s_j,1/s_j) \in \R_s \times \R_j$ as an element in $X \times X$, where $s=(s_1,s_2,\ldots)$ is an infinite sequence of positive integers, and by abuse of notation, $j$ denotes both a positive integer and a constant infinite sequence with all terms being $j$. Write $P \subseteq V \times V$ for the image of all such $p_{s,j}$'s under the map $i_X \times i_X:X \times X \to V \times V$, and $\bar{P} \subset D(V)$ for the image of $P$ under the addition. So $\bar{P} = \{[1/s_j]_s + [1/s_j]_j\}_{s,j}$.

Claim 1: $\bar{P}$ is closed in $D(V)$.

The proof goes as that of Remark~\ref{rem:closed}. In more detail, for any given plot of $V$ as in~\eqref{eq:generate}, we need to show that $B := (\rho_{s^1,\ldots,s^n})^{-1}(\bar{P})$ is closed in the domain. Let $(r_1^m,x_{1,s^1}^m,\ldots,r_n^m,x_{n,s^n}^m)$ be a sequence in $B$ with respect to $m$, which converges to $(r_1,x_{1,s^1},\ldots,r_n,x_{n,s^n})$.  This means that for each $m$, there exist an infinite sequence of positive integers $s(m)$ and a positive integer (also viewed as an infinite sequence of a constant positive integer) $j(m)$ such that $\sum_i r_i^m [x_{i,s^i}^m] = [1/s(m)_{j(m)}]_{s(m)} + [1/s(m)_{j(m)}]_{j(m)}$. One could split the set $\{1,2,\ldots,n\}$ into subsets according to the values of $x_{1,s^1},\ldots,x_{n,s^n}$, so that for $m$ large, the $x_{i,s^i}^m$'s in different subsets are mutually distinct for the same $m$. For each large $m$, there are only finitely many cases: if $s(m) \neq j(m)$, then there are two subsets, one of them has $x_{i,s^i}^m = 1/s(m)_{j(m)} \in \R_{s(m)}$ (in particular, $s^i = s(m)$ which is independent of $m$, and similarly for the rest cases which we will not mention explicitly) and the corresponding $r_i$'s adding up to $1$, and the other has $x_{k,s^k}^m = 1/s(m)_{j(m)} \in \R_{j(m)}$ and the corresponding $r_k$'s adding up to $1$, and the sum of the corresponding $r_l$'s in each of the rest subsets is $0$; if $s(m) = j(m)$, then there is a subset so that $x_{i,s^i}^m = 1/s(m)_{j(m)} \in \R_{s(m)}$ and the corresponding $r_i$'s adding up to $2$. We could pass to a subsequence so that only one fixed case occurs. So the identities keep to the limit, and hence $(r_1,x_{1,s^1},\ldots,r_n,x_{n,s^n}) \in B$.

Claim 2: $f^{-1}(\bar{P})$ is not closed in $D(V) \times D(V)$.

Corollary~\ref{cor:subspace} says that $D(X) \times D(X)$ is a topological subspace of $D(V) \times D(V)$, so it is enough to see that $f^{-1}(\bar{P}) \cap (D(X) \times D(X))$ is not closed in $D(X) \times D(X)$. Note that $([0], [0]) \notin f^{-1}(\bar{P})$, while $P \subset f^{-1}(\bar{P})$. The result follows directly from Remark~\ref{rem:notcommute} (or more precisely Example~\ref{ex:notcommute}).

These two claims then imply the result.
\end{proof}

\begin{cor}\label{prop:surprise1}\
\begin{enumerate}
\item For the diffeological vector space $V$ defined above (which is both free and projective), we have $D(V \times V) \neq D(V) \times D(V)$.
\item The projective diffeological vector space $V$ defined above is not countably-generated.
\item Write $W$ for the direct sum of $F(\R)$'s indexed over all infinite sequences of positive integers. Then $D(W)$ is not a topological vector space. In particular, $D(W \times W) \neq D(W) \times D(W)$.
\item Let $W'$ be a direct sum of $F(U)$'s indexed over a set of cardinality no less than that of the continuum, where each $U$ is open in some Euclidean space of dimension $\geq 1$. Then $D(W')$ is not a topological vector space.
\end{enumerate}
\end{cor}
\begin{proof}
(1) and (2) are straightforward from Theorem~\ref{prop:surprise2} (together with Theorem~\ref{thm:proj}).

For (3), it is enough to prove the first statement. Assume that the conclusion is not true, i.e., $D(W)$ is a topological vector space. Notice that there is a subduction $\coprod_s \R_s \to X$ of diffeological spaces, which induces a linear subduction $W \to V$. As $V$ is projective, $V$ is a smooth direct summand of $W$. Proposition~\ref{prop:retract} then implies that $V$ is a topological vector space under the $D$-topology, which conflicts with Theorem~\ref{prop:surprise2}.

For (4), we may assume that each $U$ is a ball, and the cardinality of the indexing set is equal to that of the continuum. So we have a subduction $U \to \R$ as $\dim U \geq 1$, which then induces a linear subduction $W' \to W$, where $W$ is the diffeological vector space in (3). As $W$ is projective, $W$ is then a smooth direct summand of $W'$. If $D(W')$ were a topological vector space, then so is $D(W)$ by Proposition~\ref{prop:retract}, which conflicts with (3).
\end{proof}

We summarize Proposition~\ref{prop:countable} and Corollary~\ref{prop:surprise1}(4) as follows:

\begin{thm}\label{thm:countable-vs-uncountable}
Let $Z$ be a direct sum of $F(U)$'s indexed over a set $I$, where each $U$ is open in some Euclidean space of dimension $\geq 1$. 
\begin{enumerate}
\item Then $D(Z)$ is a topological vector space if $I$ is countable, and $D(Z)$ is not a topological vector space if the cardinality of $I$ is at least that of the continuum.
\item Assuming the Continuum Hypothesis, $D(Z)$ is a topological vector space if and only if $I$ is countable.
\end{enumerate}
\end{thm}

\vspace*{10pt}

\end{document}